\documentclass[11pt]{article}

\usepackage{graphicx}
\usepackage{color,verbatim}

\usepackage{algorithm}
\usepackage{algorithmic}
\usepackage{multirow}
\usepackage{tabularx}
\renewcommand{\arraystretch}{1.0}
\newtheorem{thm}{Theorem}[section]

\newtheorem{pro}[thm]{Proposition}
\newtheorem{lem}[thm]{Lemma}
\newtheorem{cor}[thm]{Corollary}
\newtheorem{alg}[thm]{Algorithm}
\newtheorem{ass}[thm]{Assumption}
\newtheorem{defi}[thm]{Definition}
\newtheorem{example}[thm]{Example}
\newcommand{\sect}[1]{
        \par
        \stepcounter{section}
        \settowidth{\hangindent}{\large\bf\thesection.~}
        \hangafter=1
        \bigskip\bigskip\noindent
        {\large\bf\hbox{\thesection.~}#1}\par
        \nopagebreak
        \medskip
        \renewcommand{\theequation}{\thesection.\arabic{equation}}
        \setcounter{equation}{0}
        \setcounter{subsection}{0}
}
\renewcommand{\subsection}[1]{
        \stepcounter{subsection}
        \noindent
        {\bf\hbox{\thesubsection.~}#1}
        \nobreak
}
\renewcommand{\subsubsection}[1]{
        \stepcounter{subsubsection}
        \noindent
        {\bf\hbox{\thesubsubsection.~}#1}
        \nobreak
}

\newcommand{\nn}{\nonumber}
\newcommand{\be}{\begin{equation}}
\newcommand{\ee}{\end{equation}}
\newcommand{\ba}{\begin{array}}
\newcommand{\ea}{\end{array}}
\newcommand{\bea}{\begin{eqnarray}}
\newcommand{\eea}{\end{eqnarray}}
\newcommand{\bal}{\begin{alg}}
\newcommand{\eal}{\end{alg}}
\newcommand{\ble}{\begin{lem}}
\newcommand{\ele}{\end{lem}}
\newcommand{\bco}{\begin{cor}}
\newcommand{\eco}{\end{cor}}
\newcommand{\bde}{\begin{defi}}
\newcommand{\ede}{\end{defi}}
\newcommand{\bth}{\begin{thm}}
\newcommand{\eth}{\end{thm}}
\newcommand{\bpr}{\begin{pro}}
\newcommand{\epr}{\end{pro}}
\newcommand{\bas}{\begin{ass}}
\newcommand{\eas}{\end{ass}}
\newcommand{\bex}{\begin{example}}
\newcommand{\eex}{\end{example}}
\newcommand{\reff}[1]{(\ref{#1})}
\newcommand{\refa}[1]{Assumption\ \ref{#1}}
\newcommand{\refl}[1]{Lemma\ \ref{#1}}
\newcommand{\reft}[1]{Theorem\ \ref{#1}}

\newcommand{\refal}[1]{Algorithm\ \ref{#1}}

\newtheorem{prc}[thm]{Procedure}
\newcommand{\bprc}{\begin{prc}}
\newcommand{\eprc}{\end{prc}}
\def\dd{&\!\!\!\!}
\def\eop{\hfill\vbox{\hrule height0.6pt\hbox{\vrule height1.3ex
width0.6pt\hskip1.2ex\vrule width0.6pt}\hrule height0.6pt}}
\def\prf{\noindent {\sl Proof.} \rm}
\def\alglist{
\begin{list}{Step 1}
{\setlength{\leftmargin}{0.5 in}\setlength{\labelwidth}{0.7 in}}
}
\def\eli{\end{list}}
\def\bdes{\begin{description}}
\def\edes{\end{description}}

\def\na{\nabla}

\def\hf{\frac{1}{2}}

\def\st{\hbox{s.t.}}
\def\diag{\hbox{diag}\,}

\begin{document}
\pagenumbering{arabic}

\title{A primal-dual majorization-minimization method \\
for large-scale linear programs
}
\author{
Xin-Wei Liu,\thanks{Institute of Mathematics, Hebei University of Technology, Tianjin 300401, China. E-mail:
mathlxw@hebut.edu.cn. The research is supported by the Chinese NSF grants (nos. 12071108 and 11671116).} \
Yu-Hong Dai,\thanks{Academy of Mathematics and Systems Science, Chinese Academy of Sciences, Beijing 100190, China \& School of Mathematical Sciences, University of Chinese Academy of Sciences, Beijing 100049, China. This author is supported by the NSFC grants (nos. 12021001, 11991021, 11991020 and 11971372), the National Key R\&D Program of China (nos. 2021YFA 1000300 \& 2021YFA 1000301),
 the Strategic Priority Research Program of Chinese Academy of Sciences (no. XDA27000000).} \
and\
Ya-Kui Huang\thanks{Institute of Mathematics, Hebei University of Technology, Tianjin 300401, China.
This author is supported by the Chinese NSF grant (no.  11701137).} 
}
\maketitle

\noindent\underline{\hspace*{6.3in}}
\par

\vskip 10 true pt \noindent{\small{\bf Abstract.}
We present a primal-dual majorization-minimization method for solving large-scale linear programs. A smooth barrier augmented Lagrangian (SBAL) function with strict convexity for the dual linear program is derived. The majorization-minimization approach is naturally introduced to develop the smoothness and convexity of the SBAL function. Our method only depends on a factorization of the constant matrix independent of iterations and does not need any computation on step sizes, thus can be expected to be particularly appropriate for large-scale linear programs. The method shares some similar properties to the first-order methods for linear programs, but its convergence analysis is established on the differentiability and convexity of our SBAL function. The global convergence is analyzed without prior requiring either the primal or dual linear program to be feasible. Under the regular conditions, our method is proved to be globally linearly convergent, and a new iteration complexity result is given.



\noindent\underline{\hspace*{6.3in}}

\vfil\eject
}

\sect{Introduction}

We consider to solve the linear program in the dual form \bea
\min_y~-b^Ty\quad\st~A^Ty\le c, \label{prob1}\eea
where $y\in\Re^m$ is the unknown, $b\in\Re^m$, $A\in\Re^{m\times n}$ and $c\in\Re^n$ are given data. Corresponding to the dual problem \reff{prob1}, the primal linear program has the form \bea
\min_x~c^Tx\quad\st~Ax=b,~x\ge 0,\label{prob2}\eea
where $x\in\Re^n$. Problem \reff{prob2} is called the standard form of linear programming. In the literature, most of the methods and theories for linear programming are developed with the standard form (see, for example, \cite{NocWri99,SunYua06,wright97,ye}). Moreover, it is often assumed that $m<n$, $\hbox{rank}(A)=m$.

The simplex methods are the most efficient and important methods for linear programming before 1980s. These methods search the optimal solution in vertices of a polyhedral set along the boundary of the feasible region of linear programming. The initial point should be a so-called basic feasible solution corresponding to a vertex of the polyhedron which may be obtained by solving some auxiliary linear programming problem with a built-in starting point. The main computation for a new iteration point is the solution of the linear systems \bea
Bu=a,\quad B^Tv=d, \label{sec1f1}\eea
where $u\in\Re^m$ and $v\in\Re^m$ are the unknowns, $B\in\Re^{m\times m}$ is a nonsingular sub-matrix of $A$ and its one column is rotated in every iteration, $a\in\Re^m$ and $d\in\Re^m$ are some given vectors. The simplex methods are favorite since the systems in \reff{sec1f1} are thought to be easily solved.

It was discovered in \cite{KleeM72}, however, that the simplex approach could be inefficient for certain pathological problems since the number of iterations (also known as the worst-case time complexity) was exponential in the sizes of problems. In contrast, the interior-point approach initiated in 1984 by Karmarkar \cite{karmar} has been proved to be of the worst-case polynomial time complexity, a much better theoretical property than that for the simplex methods. Up to now, the best worst-case polynomial time complexity on interior-point methods is ${\cal O}(\sqrt{n}\log\frac{1}{\epsilon})$ (see, for example, \cite{wright97,ye}).

In general, interior-point methods converge to the optimal solution along a central path of the feasible polytope. The central path is usually defined by a parameter-perturbed Karush-Kuhn-Tucker (KKT) system. The system can be induced by the KKT conditions of the logarithmic-barrier problem \bea\min~c^Tx-\mu\sum_{i=1}^n\ln x_i\quad \st~Ax=b, \label{sec1f2}\eea
where $\mu>0$ is the barrier parameter, $x_i>0$ for $i=1,\ldots,n$ (that is, $x$ should be an interior-point). It is known that the well-defined central path depends on the nonempty of the set of the primal-dual interior-points \bea
{\cal F}:=\{(x,y,s)| Ax=b,~A^Ty+s=c,~x>0,~s>0\}. \nn\eea

Although there are various interior-point methods, such as the affine-scaling methods, the logarithmic-barrier methods, the potential-reduction methods, the path-following methods, etc., all these methods share some common features that distinguish them from the simplex methods.
Distinct from the simplex methods in starting from a feasible point, the interior-point methods require the initial point to be an interior-point which may not be feasible to the problem. While the simplex methods usually require a larger number of relatively inexpensive iterations, every interior-point iteration needs to solve a system with the form \bea
AS^{-1}XA^Tv=d, \label{sec1f3}\eea
where $S=\hbox{diag}(s)$ and $X=\hbox{diag}(x)$. This is generally more expensive to compute than \reff{sec1f1} but can make significant progress towards the solution. In particular, as the primal and dual iterates tend to the solutions of the primal and dual problems, some components of $x$ and $s$ can be very close to zero, which can bring about both huge and tiny values of the elements of $S^{-1}X$ and an ill-conditioned Jacobian matrix of the system \reff{sec1f3} (see \cite{NocWri99}). Some advanced methods for improving classic interior-point methods have been proposed, including the sparse matrix factorization, the Krylov subspace method and the preconditioned conjugate gradient method (see, for example, \cite{BYZ00,CMTH19,Davis06,FGZ14,Gon12a,Gon12b}).

Recently, some first-order methods for solving linear programs and linear semidefinite programming have been presented, see \cite{LMYZ,WZW} and the references therein. These methods are mainly the alternating direction augmented Lagrangian methods of multipliers (ADMM)-based methods, and can be free of solving systems \reff{sec1f1} and \reff{sec1f3}. Since the solved problems may be reformulated in different ways which result in various augmented Lagrangian function, these methods may be distinct in the augmented Lagrangian subproblems. For example, Lin et al. \cite{LMYZ} proposed their ADMM-based interior-point method based on the well-behaved homogeneous self-dual embedded linear programming model \cite{ye}, while the method in \cite{WZW} is established on using the classic augmented Lagrangian function and the projection on the cone of positive semidefinite matrices.

\subsection{Our contributions.}
We present a primal-dual majorization-minimization method on basis of solving linear programs in dual form \reff{prob1}. In our method, $y_i~(i=1,\ldots,m)$ are the primal variables, and $x_j~(j=1,\ldots,n)$ the dual variables. The method is originated from a combination of the Fiacco-McCormick logarithmic-barrier method and the Hestenes-Powell augmented Lagrangian method (see \cite{LDHmc} for more details on general nonlinear inequality-constrained optimization). A smooth barrier augmented Lagrangian (SBAL) function with strict convexity for the dual linear program is derived. Based on the smoothness and convexity of SBAL function, a majorization surrogate function is naturally designed to find the approximate minimizer of the augmented Lagrangian on primal variables, and the dual estimates are derived by a step for maximizing a minorization surrogate function of the augmented Lagrangian on dual variables. Our method can avoid the computation on the ill-conditioned Jacobian matrix like \reff{sec1f3} and does not solve some iteration-varying system \reff{sec1f1} or \reff{sec1f3} like the simplex methods and interior-point methods.

Our method initiates from the logarithmic-barrier reformulation of problem \reff{prob1}, thus can be thought of an interior-point majorization-minimization method, and shares some similar features as \cite{LMYZ}. It can also be taken as a smooth version of \cite{WZW} for linear programs, but it does not depend on any projection and computes more steps on primal iterates. Differing from the fixed-point framework for proving convergence in \cite{WZW}, based on the smoothness and convexity of our augmented Lagrangian, we can do the global convergence and prove the results on convergence rate and iteration complexity based on the well developed theories on convex optimization \cite{Nest18}.

Our proposed method only needs the factorization of the constant matrix $AA^T$, which is distinguished from the existing simplex methods and interior-point methods for linear programs necessary to solve either \reff{sec1f1} or \reff{sec1f3} varied in every iteration. Since the factorization is independent of iterations and can be done in preprocessing, our method can be implemented easily with very cheap computations, thus is especially suitable for large-scale linear programs. In addition, our method does not need any computation on step sizes, which is the other outstanding feature of our method in contrast to the existing interior-point methods for linear programs. Similar to \cite{LMYZ}, the global convergence is analyzed without prior requiring either the primal or dual linear program to be feasible. Moreover, under the regular conditions, we prove that our method can be of globally linear convergence, and a new iteration complexity result is obtained.

\subsection{Some related works.}
The augmented Lagrangian methods minimize an augmented Lagrangian function approximately and circularly with update of multipliers. The augmented Lagrangian function has been playing a very important role in the development of effective numerical methods and theories for convex and nonconvex optimization problems (see some recent references, such as \cite{BM14,CDZ15,DLS17,GKR20,GR12,GP99,HL17,HFD16,LMYZ,LiuDai18,LDHmm,WZW}). The augmented Lagrangian was initially proposed by Hestenes \cite{hesten} and Powell \cite{powell} for solving optimization problems with only equality constraints. The Hestenes-Powell augmented Lagrangian method was then generalized by Rockafellar \cite{rockaf1} to solve the optimization problems with inequality constraints. Since most of the augmented Lagrangian functions for inequality-constrained optimization depend on some kind of projection, the subproblems on the augmented Lagrangian minimization are generally solved by the first-order methods.

The majorization-minimization (MM) algorithm operates on a simpler surrogate function that majorizes the objective in minimization \cite{Lange}. Majorization can be understood to be a combination of tangency and domination. Similarly, we have the minorization-maximazation algorithm when we want to maximize an objective. The MM principle can be dated to Ortega and Rheinboldt
\cite{OR70} in 1970, where the majorization idea has been stated clearly in the context of line searches. The famed expectation-maximization (EM) principle \cite{MK08} of computational statistics is a special case of the MM principle. So far, MM methods have been developed and applied efficiently for imaging and inverse problems, computer vision problems, and so on (for example, see \cite{AIG06,cohen96,FBN06,Lange,QBP15}).

Recently, by combining the Hestenes-Powell augmented Lagrangian and the interior-point logarithmic-barrier technique (\cite{LiuSun01,LiuYua07,NocWri99,SunYua06}), the authors of \cite{LDHmc} introduce a novel barrier augmented Lagrangian function for nonlinear optimization with general inequality constraints. Distinct from the classic augmented Lagrangian function for inequality constrained optimization only first-order differentiable, the newly proposed one shares the same-order differentiability with the objective and constraint functions and is convex when the optimization is convex. In order to distinguish the new barrier augmented Lagrangian function to those proposed in \cite{GKR20,GP99}, we refer to it as the smooth barrier augmented Lagrangian (SBAL for short). For linear problems \reff{prob1} and \reff{prob2}, the SBAL functions are strictly convex and concave, respectively, with respect to the primal and dual variables. In particular, the SBAL functions are well defined without requiring either primal or dual iterates to be interior-points. These outstanding features of the SBAL functions provide natural selections for the majorization-minimization methods.

\subsection{Organization and notations.}
Our paper is organized as follows. In section 2, we describe the application of our augmented Lagrangian method in \cite{LDHmc} to the linear programs and present the associated preliminary results. The majorized functions and our primal-dual majorization-minimization method are proposed in section 3. The analysis on the global convergence and the convergence rates is done, respectively, in sections 4 and 5. We conclude our paper in the last section.

Throughout the paper, all vectors are column vectors. We use capital letters to represent matrices, and a capital letter with a subscript such as $A_i$ means the $i$th column of matrix $A$. The small letters are used to represent vectors, and a small letter with a subscript such as $s_i$ means the $i$th component of vector $s$. The capital letter $S$ means the diagonal matrix of which the components of vector $s$ are the diagonal elements. In general, we use the subscripts $k$ and $\ell$ to illustrate the letters to be related to the $k$th and $\ell$th iterations, and $i$ and $j$ the $i$th and $j$th components of a vector or the $i$th and $j$th sub-vectors of a matrix. In other cases, it should be clear from the context. To quantify the convergence of sequences, we introduce the weighted norm $\|y\|_M=\sqrt{y^TMy}$, where $y$ is a column vector, $M$ is either a positive semi-definite or positive definite symmetric matrix with the same order as $y$. The symbol $e$ is the all-one vector, for which the dimension may be varying and can be known by the context. For the symmetric positive definite matrix $B$, we use $\lambda_{\min}(B)$ and $\lambda_{\max}(B)$ to represent the minimum and maximum of eigenvalues of $B$, respectively. As usual, we use the capital letters in calligraphy to represent the index sets, $\|\cdot\|$ is the Euclidean norm, $x\circ s$ is the Hadamard product of vectors $x$ and $s$, and $x\in\Re_{++}^n$ means $x\in\Re^n$ and $x>0$ in componentwise.

\sect{The SBAL function and some preliminary results}

Recently, the authors in \cite{LDHmm,LDHmc} presented a novel barrier augmented Lagrangian function for nonlinear optimization with general inequality constraints. For problem \reff{prob1}, we reformulate it as \bea
\min_{y,s}~-b^Ty \quad\st~A^Ty+s=c,~s\ge 0, \label{sec2f1}\eea
where $s\in\Re^n$ is a slack vector. The logarithmic-barrier problem associated with \reff{sec2f1} has the form \bea
\min_{y,s}~-b^Ty-\mu\sum_{i=1}^n\ln s_i\quad\st~s-c+A^Ty=0, \label{sec2f2}\eea
where $s=(s_i)>0$, $\mu>0$ is the barrier parameter. Noting that problem \reff{sec2f2} is one with only equality constraints, we can use the Hestenes-Powell augmented Lagrangian function to reformulate it into a unconstrained optimization problem as follows, \bea
\min_{y,s}~F_{(\mu,\rho)}(y,s;x):=-\rho b^Ty-\rho\mu\sum_{i=1}^n\ln s_i+\rho x^T(s-c+A^Ty)+\hf\|s-c+A^Ty\|^2, \label{sec2f3}\eea
where $\rho>0$ is the penalty parameter which may be reduced adaptively if necessary, $x\in\Re^n$ is an estimate of the Lagrange multiplier vector.

Since $\frac{\partial^2 F_{(\mu,\rho)}(y,s;x)}{\partial s_i^2}=\frac{\rho\mu}{s_i^2}+1>0$, no matter what are $(y,s)$ and $x$, $F_{(\mu,\rho)}(y,s;x)$ is a strictly convex function with respect to $s_i$. Therefore, $F_{(\mu,\rho)}(y,s;x)$ will take the minimizer when \bea
\frac{\partial F_{(\mu,\rho)}(y,s;x)}{\partial s_i}=-\frac{\rho\mu}{s_i}+\rho x_i+(s_i-c_i+A_i^Ty)=0, \nn\eea
where $A_i\in\Re^m$ is the $i$th column vector of $A$. Equivalently, one has
 \bea s_i=\frac{1}{2}(\sqrt{(\rho x_i-c_i+A_i^Ty)^2+4\rho\mu}-(\rho x_i-c_i+A_i^Ty)). \nn\eea
Based on the observation that $s_i$ will be altered with $y$ and $x$ and is dependent on the parameters $\mu$ and $\rho$, and for simplicity of statement, we define $s=s(y,x;\mu,\rho)$ and $z=z(y,x;\mu,\rho)$ in componentwise as \bea
\dd\dd s_i(y,x;\mu,\rho)=\frac{1}{2}(\sqrt{(\rho x_i-c_i+A_i^Ty)^2+4\rho\mu}-(\rho x_i-c_i+A_i^Ty)), \label{sdf}\\
\dd\dd z_i(y,x;\mu,\rho)=\frac{1}{2}(\sqrt{(\rho x_i-c_i+A_i^Ty)^2+4\rho\mu}+(\rho x_i-c_i+A_i^Ty)), \label{zdf}\eea
where $i=1,\ldots,n$. By \reff{sdf} and \reff{zdf}, $z=s-c+A^Ty+\rho x$.
Correspondingly, the objective function $F_{(\mu,\rho)}(y,s;x)$ of the unconstrained optimization problem \reff{sec2f3} can be written as \bea
L_B(y,x;\mu,\rho)=-\rho b^Ty+\sum_{i=1}^nh_i(y,x;\mu,\rho), \label{Ldf}\eea
where $y\in\Re^m$ and $x\in\Re^n$ are the primal and dual variables of problem \reff{prob1}, $\mu>0$ and $\rho>0$ are, respectively, the barrier parameter and the penalty parameter, \bea
\dd\dd h_i(y,x;\mu,\rho)=-\rho\mu\ln s_i(y,x;\mu,\rho)+\hf z_i(y,x;\mu,\rho)^2-\frac{1}{2}\rho^2x_i^2. \label{hdf}\eea
We may write $s$ and $z$ for simplicity in the sequel when their dependence on $(y,x)$ and $(\mu,\rho)$ is clear from the context.

Similar to \cite{LDHmc}, we can prove the differentiability of the functions $s$, $z$ defined by \reff{sdf}, \reff{zdf}, and the barrier augmented Lagrangian function $L_B(y,x;\mu,\rho)$ defined by \reff{Ldf}.
\ble\label{lemzp} For given $\mu>0$ and $\rho>0$, let $L_B(y,x;\mu,\rho)$ be defined by \reff{Ldf}, $s=(s_i(y,x;\mu,\rho))\in\Re^n$ and $z=(z_i(y,x;\mu,\rho))\in\Re^n$, $S=\diag(s)$ and $Z=\diag(z)$. 

(1) Both $s$ and $z$ are differentiable with respect to $y$ and $x$, and \bea
&&\na_ys=-A(S+Z)^{-1}S, \quad \na_yz=A(S+Z)^{-1}Z, \label{20140327a}\\
&&\na_{x}s=-{\rho}(S+Z)^{-1}S, \quad
\na_{x}z={\rho}(S+Z)^{-1}Z. \label{20140327b} \eea

(2) The function $L_B(y,x;\mu,\rho)$ is twice continuously differentiable with respect to $y$, and
\bea \dd\dd \na_y L_B(y,x;\mu,\rho)=Az(y,x;\mu,\rho)-\rho b,\nn\\
     \dd\dd \na_{yy}^2L_B(y,x;\mu,\rho)=A(S+Z)^{-1}ZA^T. \nn \eea
Thus, $L_B(y,x;\mu,\rho)$ is strictly convex with respect to $y$.

(3) The function $L_B(y,x;\mu,\rho)$ is twice continuously differentiable and strictly concave with respect to $x$, and
\bea \dd\dd\na_x L_B(y,x;\mu,\rho)=\rho (s(y,x;\mu,\rho)-c+A^Ty),\nn\\
     \dd\dd \na_{xx}^2L_B(y,x;\mu,\rho)=-{\rho}^2(S+Z)^{-1}S. \nn \eea
\ele\prf (1) By \reff{sdf} and \reff{zdf}, $s-z=c-A^Ty-\rho x$ and $$s_i+z_i=\sqrt{(\rho x_i-c_i+A_i^Ty)^2+4\rho\mu}.$$ Thus, one has \bea
\dd\dd\na_y s-\na_y z=-A, \nn\\[2pt]
\dd\dd\na_y s+\na_y z=A(S+Z)^{-1}\diag(\rho x-c+A^Ty)=A(I-2(S+Z)^{-1}S). \nn\eea
Thus, by doing summation and subtraction, respectively, on both sides of the preceding equations, we have
\bea \dd\dd 2\na_y s=-2A(S+Z)^{-1}S, \nn\\[2pt]
\dd\dd -2\na_y z=-2A(I-(S+Z)^{-1}S)=-2A(S+Z)^{-1}Z. \nn \eea
Therefore, \reff{20140327a} follows immediately. The results in \reff{20140327b} can be derived in the same way by differentiating with respect to $x$.

(2) Let $h(y,x;\mu,\rho)=(h_i(y,x;\mu,\rho))\in\Re^n$. Due to \reff{hdf} and noting that $SZ=\rho\mu I$,
\bea \na_y h(y,x;\mu,\rho)=-{\rho\mu}\na_y sS^{-1}+\na_y zZ=A(S+Z)^{-1}(\rho\mu I+Z^2)=AZ. \nn
\eea Thus, $\na_yL_B(y,x;\mu,\rho)=-\rho b+\na_y h(y,x;\mu,\rho)e=Az-\rho b$. Furthermore, by (1),  \bea
\na_{yy}^2L_B(y,x;\mu,\rho)=\na_y zA^T=A(S+Z)^{-1}ZA^T. \nn \eea

(3) Note that \bea\dd\dd\na_{x} h(y,x;\mu,\rho)=-\rho\mu S^{-1}\na_{x} s+Z\na_{x} z-{\rho}^2X=\rho (Z-{\rho}X), \nn\\[2pt]
\dd\dd\na_{xx}^2 h(y,x;\mu,\rho)={\rho}(\na_x Z-\rho\na_x X), \nn\eea
and $\na_{x}L_B(y,x;\mu,\rho)=\na_{x} h(y,x;\mu,\rho)e$, $\na_{xx}^2L_B(y,x;\mu,\rho)={\rho}(\na_x z-\rho\na_x x)$.
The desired formulae in (3) can be derived immediately from the equation $s-c+A^Ty=z-\rho x$ and the results of (1). \eop

The next result gives the relation between the SBAL function and the logarithmic-barrier problem.
\bth For given $\mu>0$ and $\rho>0$, let $L_B(y,x;\mu,\rho)$ be defined by \reff{Ldf}. Then $((y^*,s^*),x^*)$ is a KKT pair of the logarithmic-barrier problem \reff{sec2f2} if and only if $s^*-c+A^Ty^*=0$ and
\bea L_B(y^*,x;\mu,\rho)\le L_B(y^*,x^*;\mu,\rho)\le L_B(y,x^*;\mu,\rho), \label{220530b}\eea
i.e., $(y^*,x^*)$ is a saddle point of the SBAL function $L_B(y,x;\mu,\rho)$. \eth\prf
Due to \refl{lemzp} (3), for any $y$ such that $c_i-A_i^Ty>0$, $L_B(y,x;\mu,\rho)$ reaches its maximum with respect to $x_i$ at $x^*_i=\frac{\mu}{c_i-A_i^Ty}$ since $\frac{\partial{L_B(y,x;\mu,\rho)}}{\partial x_i}|_{x_i=x_i^*}=0$.
If $c_i-A_i^Ty\le 0$, then $\frac{\partial L_B(y,x;\mu,\rho)}{\partial x_i}>0$, which means that $L_B(y,x;\mu,\rho)$ is strictly monotonically increasing to $\infty$ as $x_i\rightarrow \infty$.
Thus, \bea \hbox{argmax}_{x_i\in\Re} L_B(y,x;\mu,\rho)=\left\{\ba{ll}
\frac{\mu}{c_i-A_i^Ty}, & \hbox{if}\quad c_i-A_i^Ty>0; \\
\infty, & \hbox{otherwise.}\ea\right. \label{220530a}\eea

If $((y^*,s^*),x^*)$ is a KKT pair of the logarithmic-barrier problem \reff{sec2f2}, then $s^*>0$ and \bea
Ax^*=b,~s^*-c+A^Ty^*=0,~\hbox{and}~x_i^*s_i^*=\mu,~i=1,\ldots,n. \nn\eea
Thus, $s_i^*=c_i-A_i^Ty^*>0$ and $x_i^*=\frac{\mu}{c_i-A_i^Ty^*},~i=1,\ldots,n$.
Therefore, by \reff{220530a}, \bea L_B(y^*,x^*;\mu,\rho)=-\rho b^Ty^*-\rho\mu\sum_{i=1}^n\ln (c_i-A_i^Ty^*)\ge L_B(y^*,x;\mu,\rho). \nn\eea
Furthermore, the condition $x_i^*s_i^*=\mu$ implies $z_i(y^*,x^*;\mu,\rho)-\rho x_i^*=0.$ Thus, $Az(y^*,x^*;\mu,\rho)=\rho b.$
It follows from \refl{lemzp} (2), $y^*$ is the minimizer of $L_B(y,x^*;\mu,\rho)$. That is, the right-hand-side inequality in \reff{220530b} holds.

In reverse, if $(y^*,x^*)$ satisfies \reff{220530b}, then $y^*$ is a minimizer of $L_B(y,x^*;\mu,\rho)$ and $x^*$ is a maximizer of $L_B(y^*,x;\mu,\rho)$. Thus, due to \refl{lemzp} (2) and (3), one has \bea
Az(y^*,x^*;\mu,\rho)=\rho b,\quad s(y^*,x^*;\mu,\rho)-c+A^Ty^*=0. \nn\eea
The second equation further implies $z(y^*,x^*;\mu,\rho)-\rho x^*=0$ and $x_i^*(c_i-A_i^Ty^*)=\mu,~i=1,\ldots,n$.
Let $s^*=s(y^*,x^*;\mu,\rho)$. Then $s^*=c-A^Ty^*$, and $((y^*,s^*),x^*)$ is a KKT pair of the logarithmic-barrier problem \reff{sec2f2}.
\eop

The following result shows that, under suitable conditions, a minimizer of problem \reff{prob1} is an approximate minimizer of the SBAL function.
\bth\label{a} Let $y^*$ be a minimizer of the problem \reff{prob1} and $x^*$ is the associated Lagrange multiplier vector. If the Slater constraint qualification holds, then for $\mu>0$ sufficiently small and for $\rho>0$, there exists a neighborhood of $x^*$ such that $y^*$ is a $\sqrt{\rho\mu}$-approximate strict global minimizer of the augmented Lagrangian $L_B(y,x;\mu,\rho)$ (that is, there a scalar $\delta>0$ such that $\|\na_y L_B(y^*,x;\mu,\rho)\|\le\delta\sqrt{\rho\mu}$).
\eth\prf Under the conditions of the theorem, $x^*$ is a KKT point of problem \reff{prob1}. Thus, \bea
Ax^*=b,~A^Ty^*\le c,~x^*\ge 0,\ (x^*)^T(c-A^Ty^*)=0. \label{20210603a}\eea
Let $z^*_i=z_i(y^*,x^*;\mu,\rho)$. Note that $x_i^*(c_i-A_i^Ty^*)=0$ for $i=1,\ldots,n$. Then \bea
z_i^*=\left\{\ba{ll}
\hf (\sqrt{(\rho x_i^*)^2+4\rho\mu}+\rho x_i^*), & \hbox{if}\ c_i-A_i^Ty^*=0,~x_i^*>0; \\[5pt]
\hf (\sqrt{(c_i-A_i^Ty^*)^2+4\rho\mu}-(c_i-A_i^Ty^*)), & \hbox{if}\ c_i-A_i^Ty^*>0,~x_i^*=0;\\[5pt]
\sqrt{\rho\mu}, &\hbox{otherwise.} \ea\right.\nn\eea
Since $\sqrt{(\rho x_i^*)^2+4\rho\mu}\le\rho x_i^*+2\sqrt{\rho\mu}$ and $\sqrt{(c_i-A_i^Ty^*)^2+4\rho\mu}\le (c_i-A_i^Ty^*)+2\sqrt{\rho\mu}$, one has \bea
\rho x^*\le z^*\le\rho x^*+\sqrt{\rho\mu},\quad \|z^*-\rho x^*\|_{\infty}\le \sqrt{\rho\mu}. \label{20210603b}\eea

We will prove the result by showing $\|\na_y L_B(y^*,x^*;\mu,\rho)\|\le\delta\sqrt{\rho\mu}$ for some scalar $\delta$ and $\na^2_{yy} L_B(y^*,x^*;\mu,\rho)$ is positive definite for $\rho>0$. By using \refl{lemzp} (2), and \reff{20210603a}, \reff{20210603b}, we have \bea
\|\na_y L_B(y^*,x^*;\mu,\rho)\|=\|Az^*-\rho b\|=\|A(z^*-\rho x^*)\|\le\sqrt{\rho\mu}\|A\|_1, \nn \eea
which verifies the first part of the result.

Now we prove the second part of the result by showing that $d^T\na^2_{yy} L_B(y^*,x^*;\mu,\rho)d>0$ for all nonzero $d\in\Re^n$ and $\rho>0$. Let $s_i^*=s_i(y^*,x^*;\mu,\rho)$. Then \bea
\frac{z_i^*}{s_i^*+z_i^*}=\left\{\ba{ll}
\hf (1+\frac{\rho x_i^*}{\sqrt{(\rho x_i^*)^2+4\rho\mu}}), & \hbox{if}\ c_i-A_i^Ty^*=0,~x_i^*>0; \\[5pt]
\hf (1-\frac{(c_i-A_i^Ty^*)}{\sqrt{(c_i-A_i^Ty^*)^2+4\rho\mu}}), & \hbox{if}\ c_i-A_i^Ty^*>0,~x_i^*=0;\\[5pt]
\hf, &\hbox{otherwise.} \ea\right.\nn\eea
Therefore, by \refl{lemzp} (2), \bea
\dd\dd\na^2_{yy} L_B(y^*,x^*;\mu,\rho) \nn\\
\dd\dd=\sum_{i=1}^n \frac{z_i^*}{s_i^*+z_i^*}A_iA_i^T\quad (A_i~\hbox{is the}~i\hbox{th column of}~A) \nn\\
\dd\dd=\hf(\sum_{i\in {\cal I}_1}(1+\frac{\rho x_i^*}{\sqrt{(\rho x_i^*)^2+4\rho\mu}})A_iA_i^T+\sum_{i\in {\cal I}_2}A_iA_i^T) \nn\\
\dd\dd\quad+\hf\sum_{i\in {\cal I}_3}(1-\frac{(c_i-A_i^Ty^*)}{\sqrt{(c_i-A_i^Ty^*)^2+4\rho\mu}})A_iA_i^T \nn\\
\dd\dd\ge\hf(1-\max\{\frac{(c_i-A_i^Ty^*)}{\sqrt{(c_i-A_i^Ty^*)^2+4\rho\mu}},~i=1,\ldots,n\})AA^T, \eea
where ${\cal I}_1=\{i|c_i-A_i^Ty^*=0,~x_i^*>0\}$, ${\cal I}_2=\{i|c_i-A_i^Ty^*=0,~x_i^*=0\}$, ${\cal I}_3=\{i|c_i-A_i^Ty^*>0,~x_i^*=0\}$. The result follows easily because of the positive definiteness of $AA^T$. \eop

Based on the newly proposed barrier augmented Lagrangian function, \cite{LDHmc} presented a novel augmented Lagrangian method of multipliers for optimization with general inequality constraints. The method alternately updates the primal and dual iterates by \bea
\dd\dd y_{k+1}=\hbox{argmin}_y L_B(y,x_k;\mu_k,\rho_k), \label{subp1}\\
\dd\dd x_{k+1}=\frac{1}{\rho_k} z(y_{k+1},x_k;\mu_k,\rho_k). \label{subp2}\eea
The update of parameters $\mu_{k+1}$ and $\rho_{k+1}$ depend on the residual $\|s(y_{k+1},x_{k+1};\mu_k,\rho_k)-c+A^Ty_{k+1}\|$ and the norm $\|x_{k+1}\|$ of dual multiplier vector. 

To end this section, we show some monotone properties of our defined functions $L_B(y,x;\mu,\rho)$, $s_i(y,x;\mu,\rho)$ and $z_i(y,x;\mu,\rho)$ with respect to the parameters.
\ble\label{lem24an} Denote $L_B(y,x;\mu,\rho)=\rho\phi(y,x;\mu,\rho)+\hf R^2(y,x;\mu,\rho)$, where \bea
\dd\dd \phi(y,x;\mu,\rho)=-b^Ty-\mu\sum_{i=1}^n\ln s_i(y,x;\mu,\rho)+x^T(s(y,x;\mu,\rho)-c+A^Ty), \nn\\
\dd\dd R(y,x;\mu,\rho)=\|s(y,x;\mu,\rho)-c+A^Ty\|. \nn\eea
Let $\hat y_{k+1}=\hbox{\rm argmin}_y L_B(y,x_k;\mu_k,\hat\rho_k)$ and $\tilde y_{k+1}=\hbox{\rm argmin}_y L_B(y,x_k;\mu_k,\tilde\rho_k)$ be attained. If $\hat\rho_k>\tilde\rho_k$, then \bea \phi(\hat y_{k+1},x_k;\mu_k,\hat\rho_k)<\phi(\tilde y_{k+1},x_k;\mu_k,\tilde\rho_k),~
R(\hat y_{k+1},x_k;\mu_k,\hat\rho_k)>R(\tilde y_{k+1},x_k;\mu_k,\tilde\rho_k). \nn
\eea\ele\prf Let $\hat s_{k+1}=s(\hat y_{k+1},x_k;\mu_k,\hat\rho_k)$ and $\tilde s_{k+1}=s(\tilde y_{k+1},x_k;\mu_k,\tilde\rho_k)$. Then, by \reff{sec2f3}, \bea (\hat y_{k+1},\hat s_{k+1})=\hbox{\rm argmin}_{y,s}F_{(\mu_k,\hat\rho_k)}(y,s;x_k),~(\tilde y_{k+1},\tilde s_{k+1})=\hbox{\rm argmin}_{y,s}F_{(\mu_k,\tilde\rho_k)}(y,s;x_k). \nn \eea
Thus, if we denote $\psi_{\mu}(y,s;x)=-b^Ty-\mu\sum_{i=1}^n s_i+x^T(s-c+A^Ty)$ and $W(y,s;x)=\|s-c+A^Ty\|$, then
$F_{(\mu,\rho)}(y,s;x)=\rho\psi_{\mu}(y,s;x)+\hf W^2(y,s;x)$, and
\bea \dd\dd \phi(\hat y_{k+1},x_k;\mu_k,\hat\rho_k)=\psi_{\mu_k}(\hat y_{k+1},\hat s_{k+1};x_k),~\phi(\tilde y_{k+1},x_k;\mu_k,\tilde\rho_k)=\psi_{\mu_k}(\tilde y_{k+1},\tilde s_{k+1};x_k), \nn\\
\dd\dd R(\hat y_{k+1},x_k;\mu_k,\hat\rho_k)=W(\hat y_{k+1},\hat s_{k+1};x_k),~R(\tilde y_{k+1},x_k;\mu_k,\tilde\rho_k)=W(\tilde y_{k+1},\tilde s_{k+1};x_k). \nn \eea
Moreover, \bea \dd\dd F_{(\mu_k,\hat\rho_k)}(\hat y_{k+1},\hat s_{k+1};x_k)< F_{(\mu_k,\hat\rho_k)}(\tilde y_{k+1},\tilde s_{k+1};x_k), \nn\\
\dd\dd F_{(\mu_k,\tilde\rho_k)}(\tilde y_{k+1},\tilde s_{k+1};x_k)< F_{(\mu_k,\tilde\rho_k)}(\hat y_{k+1},\hat s_{k+1};x_k). \nn\eea
It follows that \bea
\dd\dd F_{(\mu_k,\hat\rho_k)}(\tilde y_{k+1},\tilde s_{k+1};x_k)-F_{(\mu_k,\hat\rho_k)}(\hat y_{k+1},\hat s_{k+1};x_k) \nn\\
\dd\dd +F_{(\mu_k,\tilde\rho_k)}(\hat y_{k+1},\hat s_{k+1};x_k)-F_{(\mu_k,\tilde\rho_k)}(\tilde y_{k+1},\tilde s_{k+1};x_k) \nn\\
\dd\dd =(\hat\rho_k-\tilde\rho_k)(\psi_{\mu_k}(\tilde y_{k+1},\tilde s_{k+1};x_k)-\psi_{\mu_k}(\hat y_{k+1},\hat s_{k+1};x_k))>0. \nn\eea
That is, $\phi(\tilde y_{k+1},x_k;\mu_k,\tilde\rho_k)=\psi_{\mu_k}(\tilde y_{k+1},\tilde s_{k+1};x_k)>\psi_{\mu_k}(\hat y_{k+1},\hat s_{k+1};x_k)=\phi(\hat y_{k+1},x_k;\mu_k,\hat\rho_k)$. Therefore, \bea
\dd\dd \hf R^2(\hat y_{k+1},x_k;\mu_k,\hat\rho_k)=\hf W^2(\hat y_{k+1},\hat s_{k+1};x_k) \nn\\
\dd\dd >\hf W^2(\hat y_{k+1},\hat s_{k+1};x_k)+\tilde\rho_k (\psi_{\mu_k}(\hat y_{k+1},\hat s_{k+1};x_k)-\psi_{\mu_k}(\tilde y_{k+1},\tilde s_{k+1};x_k)) \nn\\
\dd\dd >\hf W^2(\tilde y_{k+1},\tilde s_{k+1};x_k)=\hf R^2(\tilde y_{k+1},x_k;\mu_k,\tilde\rho_k), \nn\eea
which completes our proof. \eop


\ble\label{lem24n} For given parameters $\mu>0$ and $\rho>0$, the following results are true.

(1) Both $s_i(y,x;\mu,\rho)$ and $z_i(y,x;\mu,\rho)$ are monotonically increasing with respect to $\mu$.

(2) If $(s_i(y,x;\mu,\rho)-c_i+A_i^Ty)\ne 0$, then $(s_i(y,x;\mu,\rho)-c_i+A_i^Ty)^2$ will be decreasing as $\rho$ is decreasing.

(3) If $\|s(y,x;\mu,\rho)-c+A^Ty\|\ne 0$, then the function $\frac{1}{\rho}L_B(y,x;\mu,\rho)$ is monotonically decreasing with respect to $\rho$.

(4) The function $\frac{1}{\rho}L_B(y,x;\mu,\rho)$ is strictly convex with respect to $\rho$.
\ele\prf It should be noted that all related functions are differentiable with respect to $\mu$ and $\rho$. In addition, due to \bea
\dd\dd\frac{\partial s_i(y,x;\mu,\rho)}{\partial\mu}=\frac{\partial z_i(y,x;\mu,\rho)}{\partial\mu}=\frac{\rho}{s_i+z_i}>0, \\[2pt]
\dd\dd\frac{\partial ((s_i(y,x;\mu,\rho)-c_i+A_i^Ty)^2)}{\partial\rho}=\frac{2}{\rho}\frac{s_i}{s_i+z_i}(s_i-c_i+A_i^Ty)^2>0, \\[2pt]
\dd\dd\frac{\partial \frac{1}{\rho}L_B(y,x;\mu,\rho)}{\partial\rho}=-\frac{1}{2}\frac{1}{\rho^2}\|s-c+A^Ty\|^2<0, \\[2pt]
\dd\dd\frac{\partial^2 (\frac{1}{\rho}L_B(y,x;\mu,\rho))}{\partial\rho^2}=\frac{1}{\rho^3}(s-c+A^Ty)^T(S+Z)^{-1}Z(s-c+A^Ty), \eea
the desired results are obtained immediately. \eop

\sect{Our primal-dual majorization-minimization method}

Our method in this paper focuses on how to solve the subproblem \reff{subp1} efficiently. Noting the strict convexity of the SBAL function $L_B(y,x;\mu,\rho)$ with respect to $y$ and the special structure of the Hessian matrix $\na_{yy}^2L(y,x;\mu,\rho)$, the introduction of the majorization-minimization method is a natural selection. In particular, we will see that the dual update is precisely a step which can be derived by the minorization-maximization.

Let $(y_k,x_k)$ be the current iteration point, $\mu_k>0$ and $\rho_k>0$ are the current values of the parameters. For any given $x\in\Re^n$, we consider the quadratic surrogate function $Q_k(\cdot,x): \Re^m\to\Re$, \bea
Q_k(y,x)=\dd\dd L_B(y_k,x;\mu_k,\rho_k)+(Az(y_k,x;\mu_k,\rho_k)-\rho_k b)^T(y-y_k) \nn\\
       \dd\dd +\hf(y-y_k)^TAA^T(y-y_k), \label{Qf}\eea
which is an approximate function of the objective in \reff{subp1} and majorizes the objective function with respect to $y$.
\ble\label{lem21} For any given $x=\hat x$ and the parameters $\mu_k>0$ and $\rho_k>0$, there holds $Q_k(y_k,\hat x)=L_B(y_k,\hat x;\mu_k,\rho_k)$ and $L_B(y,\hat x;\mu_k,\rho_k)\le Q_k(y,\hat x)$ for all $y\in\Re^m$. \ele\prf
The equation $Q_k(y_k,\hat x)=L_B(y_k,\hat x;\mu_k,\rho_k)$ is obtained from \reff{Qf}.

By Taylor's theorem with remainder, \bea
\dd\dd L_B(y,\hat x;\mu_k,\rho_k)=L_B(y_k,\hat x;\mu_k,\rho_k)+\na_y L_B(y_k,\hat x;\mu_k,\rho_k)^T(y-y_k) \nn\\
\dd\dd \quad+\int_0^1\left(\na_yL_B(y_k+\tau (y-y_k),\hat x;\mu_k,\rho_k)-\na_y L_B(y_k,\hat x;\mu_k,\rho_k)\right)^T(y-y_k)d\tau. \label{lem21f1} \eea
Due to \refl{lemzp} (2), one has \bea
\dd\dd\na_yL_B(y_k+\tau (y-y_k),\hat x;\mu_k,\rho_k)-\na_y L_B(y_k,\hat x;\mu_k,\rho_k) \nn\\
\dd\dd=\int_{0}^1\tau\na_{yy}^2L_B(y_k+\alpha\tau (y-y_k),\hat x;\mu_k,\rho_k)(y-y_k)d\alpha \nn\\
\dd\dd=\int_{0}^1\tau A(\hat S_k+\hat Z_k)^{-1}\hat Z_kA^T(y-y_k)d\alpha \nn\\
\dd\dd=\tau AA^T(y-y_k)-\int_{0}^1\tau A(\hat S_k+\hat Z_k)^{-1}\hat S_kA^T(y-y_k)d\alpha, \nn \eea
where $\hat S_k=\diag(s(y_k+\tau (y-y_k),\hat x;\mu_k,\rho_k))$ and $\hat Z_k=\diag(z(y_k+\tau (y-y_k),\hat x;\mu_k,\rho_k))$.
Noting $$\int_0^1\int_{0}^1\tau (y-y_k)^TA(\hat S_k+\hat Z_k)^{-1}\hat S_kA^T(y-y_k)d\alpha d\tau\ge 0,$$
the inequality $L_B(y,\hat x;\mu_k,\rho_k)\le Q_k(y,\hat x)$ follows from \refl{lemzp} (2) and \reff{lem21f1} immediately. \eop

In a similar way, if for given $y\in\Re^m$ and the parameters $\mu_k>0$ and $\rho_k>0$, we define $P_k(y,\cdot): \Re^n\to\Re$ be the function \bea
P_k(y,x)=\dd\dd L_B(y,x_k;\mu_k,\rho_k)+\rho_k(s(y,x_k;\mu_k,\rho_k)-c+A^Ty)^T(x-x_k) \nn\\
       \dd\dd -\frac{1}{2}{\rho_k}^2(x-x_k)^T(x-x_k), \label{Df}\eea
then $P_k(y,x_k)=L_B(y,x_k;\mu_k,\rho_k)$ and $L_B(y,x;\mu_k,\rho_k)\ge P_k(y,x)$ for all $x\in\Re^n$. That is, $P_k(y,x)$ is an approximate surrogate function of the objective in optimization \bea
\max_x L_B(y,x;\mu_k,\rho_k) \nn\eea
and minorizes the objective function with respect to $x$ (i.e., majorizes the negative objective function).

By the strict convexity of $Q_k(\cdot,x)$ and the strict concavity of $P_k(y,\cdot)$, there are a unique minimizer of $Q_k(y,\hat x)$ and a unique maximizer of $P_k(\hat y,x)$, where $\hat x\in\Re^n$ and $\hat y\in\Re^m$ are any given vectors.
\ble\label{lem22} Given $\mu_k>0$ and $\rho_k>0$. Let $Q_k(\cdot,x): \Re^m\to\Re$ and $P_k(y,\cdot): \Re^n\to\Re$ be functions defined by \reff{Qf} and \reff{Df}, respectively.

(1) For any given $\hat x$, $Q_k(y,\hat x)$ has a unique minimizer $y^*_k$. Moreover, $y_k^*$ satisfies the equation \bea
AA^T(y-y_k)=-(Az(y_k,\hat x;\mu_k,\rho_k)-\rho_k b).\eea

(2) For any given $\hat y$, $P_k(\hat y,x)$ has a unique maximizer $x_k^*$, and \bea
 x_k^*=x_k+\frac{1}{\rho_k} (s(\hat y,x_k;\mu_k,\rho_k)-c+A^T\hat y). \label{220601a}\eea

(3) For any given $\hat x$ and $\hat y$, one has \bea
\dd\dd L_B(y_k^*,\hat x;\mu_k,\rho_k)-L_B(y_k,\hat x;\mu_k,\rho_k)\le-\frac{1}{2}\|Az(y_k,\hat x;\mu_k,\rho_k)-\rho_k b\|_{(AA^T)^{-1}}^2, \label{lem22f3}\\
\dd\dd L_B(\hat y,x_k^*;\mu_k,\rho_k)-L_B(\hat y,x_k;\mu_k,\rho_k)\ge\frac{1}{2}\|s(\hat y,x_k;\mu_k,\rho_k)-c+A^T\hat y\|^2. \eea
\ele\prf Since  \bea
\dd\dd\na_yQ_k(y,\hat x)=AA^T(y-y_k)+(Az(y_k,\hat x;\mu_k,\rho_k)-{\rho_k} b), \nn\\
\dd\dd\na_xP_k(\hat y,x)=-{\rho_k}^2(x-x_k)+\rho_k(s(\hat y,x_k;\mu_k,\rho_k)-c+A^T\hat y), \nn\eea
and noting the strict convexity of $Q_k(y,\hat x)$ with respect to $y$, and the strict concavity of $P_k(\hat y,x)$ with respect to $x$, the results (1) and (2) are obtained immediately from the optimality conditions of general unconstrained optimization (see \cite{NocWri99,SunYua06}).

By the preceding results, one has \bea
\dd\dd Q_k(y_k^*,\hat x)=Q_k(y_k,\hat x)-\frac{1}{2}\|Az(y_k,\hat x;\mu_k,\rho_k)-\rho_kb\|_{(AA^T)^{-1}}^2, \nn\\
\dd\dd P_k(\hat y,x_k^*)=P_k(\hat y,x_k)+\hf\|s(\hat y,x_k;\mu_k,\rho_k)-c+A^T\hat y\|^2. \nn\eea
Due to \refl{lem21}, there hold \bea
\dd\dd L_B(y_k^*,\hat x;\mu_k,\rho_k)-L_B(y_k,\hat x;\mu_k,\rho_k)\le Q_k(y_k^*,\hat x)-Q_k(y_k,\hat x), \nn\\
\dd\dd L_B(\hat y,x_k^*;\mu_k,\rho_k)-L_B(\hat y,x_k;\mu_k,\rho_k)\ge P_k(\hat y,x_k^*)-P_k(\hat y,x_k), \nn\eea
which complete our proof. \eop

Because of \reff{sdf} and \reff{zdf}, \reff{220601a} is equivalent to $x_k^*=\frac{1}{\rho_k} z(\hat y,x_k;\mu_k,\rho_k)$, which is consistent with \reff{subp2}. This fact shows that the dual update $x_{k+1}$ in \reff{subp2} can be obtained from maximizing the minorized function $P_k(y_{k+1},x)$. In the following, we describe our algorithm for linear programming.

\noindent\underline{\hspace*{6.3in}}\\[-10pt]
\bal\label{alg1}(A primal-dual majorization-minimization method for problem \reff{prob1})
{\small \alglist
\item[{\bf Step}] {\bf 0}. Given $(y_0,x_0)\in\Re^{m}\times\Re^n$, $\mu_0>0$, $\rho_0>0$, $\delta>0$, $\gamma\in(0,1)$, $\epsilon>0$.
Set $k:=0$.

\item[{\bf Step}] {\bf 1}. Approximately minimize $L_B(y,x_k;\mu_k,\rho_k)$ by the majorization-minimization method starting from $y_k$.

Set $\hat y_{0}=y_k$, $\hat\rho_{0}=\rho_k$, $\ell:=0$.

{\bf Step 1.1}. Solve the equation \bea
 AA^T(y-\hat y_{\ell})=-(Az(\hat y_{\ell},x_k;\mu_{k},\hat\rho_{\ell})-\hat\rho_{\ell}b) \label{yiter}\eea
\hspace{0.5cm} to obtain the solution $\hat y_{\ell+1}$. Evaluate \bea
E_{k+1}^{primal}=\|A z(\hat y_{\ell+1},x_{k};\mu_k,\hat\rho_{\ell})-\hat\rho_{\ell}b\|. \nn\eea
\hspace{0.5cm} If $E_{k+1}^{primal}>\mu_{k}$,
set $\hat\rho_{\ell+1}=\hat\rho_{\ell}$, $\ell:=\ell+1$ and repeat Step 1.1. Otherwise, compute \bea
E_{k+1}^{dual}=\|s(\hat y_{\ell+1},x_{k};\mu_{k},\hat\rho_{\ell})-c+A^T\hat y_{\ell+1}\|. \nn\eea
\hspace{0.5cm} If $E_{k+1}^{dual}>\max\{\hat\rho_{\ell},\mu_k\}$, set $\hat\rho_{\ell+1}\ge 0.5\hat\rho_{\ell}$,  $\ell:=\ell+1$ and repeat Step 1.1; else set
$y_{k+1}=\hat y_{\ell+1}$,

\hspace{0.5cm} $\rho_{k+1}=\hat\rho_{\ell}$, end.

\item[{\bf Step}] {\bf 2}. Update $x_k$ to \bea
x_{k+1}=x_{k}+\frac{1}{\rho_{k+1}}(s(y_{k+1},x_{k};\mu_{k},\rho_{k+1})-c+A^Ty_{k+1}). \label{xiter}\eea

\item[{\bf Step}] {\bf 3}. If $\mu_k<\epsilon$, stop the algorithm. Otherwise,
set $\mu_{k+1}\le\gamma\mu_k$, $\rho_{k+1}=\min\{\rho_{k+1}, \frac{\delta}{\|x_{k+1}\|_{\infty}}\}$, $k:=k+1$. End (while)


\eli}
\eal
\noindent\underline{\hspace*{6.3in}}

The initial point for our algorithm can be arbitrary, which is different from both the simplex methods and the interior-point methods starting from either a feasible point or an interior-point. Theoretically, since the augmented Lagrangian function is an exact penalty function, we can always select the initial penalty parameter $\rho_0$ sufficiently small such that, under desirable conditions, $E_{k+1}^{dual}$ is sufficiently small. The initial barrier parameter $\mu_0$ can be selected to be small without affecting the well-definedness of the algorithm, but it may impact the strict convexity of the SBAL function and bring about more iterations for solving the subproblem \reff{subp1}.

The Step 1 is the core and the main computation of our algorithm. For fixed $x_k$, $\mu_k$ and $\rho_k$, we attempt to find a new estimate $y_{k+1}$, which is an approximate minimizer of the SBAL function $L_B(y,x_k;\mu_k,\rho_k)$ with respect to $y$. The main computation is in solving the system \reff{yiter}, which depends on the decomposition of $AA^T$. Since $AA^T$ is independent of the iteration, its decomposition can be fulfilled in preprocessing. If $L_B(y,x_k;\mu_k,\rho_k)$ is lower bounded, then the Step 1 will terminate in a finite number of iterations.

By Step 2 of \refal{alg1}, we have $x_{k+1}=\frac{1}{\rho_{k+1}} z(y_{k+1},x_{k};\mu_{k},\rho_{k+1})$, thus $x_{k+1}>0$ for all $k\ge 0$. Due to \refl{lem22} (3) and the strict concavity, one has \bea \|s(y_{k+1},x_{k+1};\mu_{k},\rho_{k+1})-c+A^Ty_{k+1}\|<\|s(y_{k+1},x_{k};\mu_{k},\rho_{k+1})-c+A^Ty_{k+1}\|.
\label{220612c}\eea
Due to the Step 3, $\mu_k\to 0$ as $k\to\infty$, $\rho_{k+1}\|x_{k+1}\|_{\infty}\le\delta$ for all $k>0$.

\sect{Global convergence}

We analyze the convergence of \refal{alg1} in this section. Firstly, we prove that, if the original problem has a minimizer, then the Step 1 will always terminate in a finite number of iterations and $\{y_k\}$ will be obtained. After that, we prove that, without prior requiring either the primal or the dual linear problem to be feasible, our algorithm can recognize the KKT point of problem \reff{prob1}, or illustrate that either its dual problem \reff{prob2} is unbounded as problem \reff{prob1} is feasible, or a minimizer with lease violations of constraints is found as problem \reff{prob1} is infeasible.

\ble\label{le41n} If problem \reff{prob1} has a solution, then for any given $x_k\in\Re^n_{++}$ and any given parameters $\mu_k>0$ and $\rho_k>0$, the SBAL function $L_B(y,x_k;\mu_k,\rho_k)$ is lower bounded from $-\infty$, and the Step 1 will terminate in a finite number of iterations. \ele\prf If problem \reff{prob1} has a solution, then the logarithmic-barrier problem \reff{sec2f2} is feasible when the original problem is strictly feasible (that is, the Slater constraint qualification holds), otherwise problem \reff{sec2f2} is infeasible. Correspondingly, the objective
$-b^Ty-\mu\sum_{i=1}^n\ln s_i$ of problem \reff{sec2f2} either takes its minimizer at an interior-point of problem \reff{prob1} (in this case the minimizer is attained) or is $+\infty$. It is noted that $A_i^Ty\rightarrow-\infty$ for any $i=1,\ldots,n$ if and only if $A_i^Ty<c_i$, the strict feasibility of the corresponding constraint of problem \reff{prob1}. The preceding result shows that no matter when $y$ is such that $A_i^Ty\to-\infty$ for any $i=1,\ldots,n$, the minimizer of $-b^Ty-\mu\sum_{i=1}^n\ln s_i$ with $s_i=\max\{c_i-A_i^Ty,0\}$ will be lower bounded away from $-\infty$.

If $L_B(y,x_k;\mu_k,\rho_k)$ is not lower bounded, then $L_B(y,x_k;\mu_k,\rho_k)\to-\infty$ as $A_i^Ty\rightarrow-\infty$ for some $i=1,\ldots,n$. Let ${\cal I}(y)=\{i|A_i^Ty\to-\infty\}$. Since \bea
\dd\dd L_B(y,x_k;\mu_k,\rho_k)\ge-b^Ty-\mu_k\sum_{i=1}^n\ln s_i(y,x_k;\mu_k,\rho_k) \nn\\
\dd\dd=-b^Ty-\mu_k\sum_{i\in{\cal I}(y)}\ln (c_i-A_i^Ty)-\mu_k\sum_{i\in{\cal I}(y)}\ln\frac{s_i(y,x_k;\mu_k,\rho_k)}{c_i-A_i^Ty}-\mu_k\sum_{i\not\in{\cal I}(y)}\ln s_i(y,x_k;\mu_k,\rho_k) \nn\\
\dd\dd>-\infty, \nn\eea
it shows that $L_B(y,x_k;\mu_k,\rho_k)$ is lower bounded away from $-\infty$.

Now we prove that for any fixed $\hat\rho_{\ell}$, if the Step 1 of \refal{alg1} does not terminate finitely, then $E_{k+1}^{primal}\to 0$ as $\ell\to\infty$. By \refl{lem22}, $\{L_B(\hat y_{\ell},x_k;\mu_k,\hat\rho_{\ell})\}$ is monotonically non-increasing as $\ell\to\infty$. Thus either there is a finite limit for the sequence $\{L_B(\hat y_{\ell},x_k;\mu_k,\hat\rho_{\ell})\}$ or the whole sequence tends to $-\infty$. Since $L_B(y,x_k;\mu_k,\hat\rho_{\ell})$ is bounded below, due to \reff{lem22f3}, one has \bea
\lim_{\ell\to\infty} \|Az(\hat y_{\ell},x_k;\mu_k,\hat\rho_{\ell})-\hat\rho_{\ell} b\|_{(AA^T)^{-1}}=0, \eea
which shows that the condition $E_{k+1}^{primal}\le\mu_k$ will be satisfied in a finite number of iterations.

Since problem \reff{prob1} is supposed to be feasible, for every $s>0$ one has $s-c+A^Ty>0$. It follows from \refl{lem24an} that there is a scalar $\rho_{k+1}>0$ such that for given $\mu_k>0$ and for all $\hat\rho_{\ell}\le\rho_{k+1}$, $E_{k+1}^{dual}\le\mu_k$ as $\ell$ is large enough. Thus, the Step 1 will terminate in a finite number of iterations.
\eop

The next result shows that, if the Step 1 does not terminate finitely, then either problem \reff{prob1} is unbounded or a point with least constraint violations will be found.
\ble\label{le42n} For given $x_k\in\Re_{++}^n$ and parameters $\mu_k>0$ and $\rho_k>0$, if the Step 1 of \refal{alg1} does not terminate finitely and an infinite sequence $\{\hat y_{\ell}\}$ is generated, then either problem \reff{prob1} is unbounded or any cluster point of $\{\hat y_{\ell}\}$ is an infeasible stationary point $y^*$ satisfying \bea
A\max\{A^Ty^*-c,0\}=0. \eea \label{s4f2a}
The point $y^*$ is also a solution for minimizing the $\ell_2$-norm of constraint violations of problem \reff{prob1}, and shows that problem \reff{prob2} is unbounded.\ele\prf If that the Step 1 of \refal{alg1} does not terminate finitely is resulted from $E_{k+1}^{primal}$ not being small enough for given $\rho_k$, then $\{\hat y_{\ell}\}$ is unbounded and $L_B(\hat y_{\ell},x_k;\mu_k,\hat\rho_{\ell})\to-\infty$ as $\ell\to\infty$, which by the arguments in the proof of \refl{le41n} implies that problem \reff{prob1} is feasible and unbounded.

Now we consider the case that $\{\hat y_{\ell}\}$ is bounded for given $\rho_k$. Suppose that ${\cal L}$ is a subset of indices such that $\hat y_{\ell}\to{\hat y_{\ell}}^*$ as $\ell\in{\cal L}$ and $\ell\to\infty$ for given $\rho_k$. Then  \bea
Az({\hat y_{\ell}}^*,x_k;\mu_k,\rho_k)-\rho_k b=0. \label{220725a}\eea
Due to $z({\hat y_{\ell}}^*,x_k;\mu_k,\rho_k)>0$, \reff{220725a} shows that problem \reff{prob2} is feasible. Furthermore, considering the fact that the Step 1 of \refal{alg1} does not terminate finitely, one has $\hat\rho_{\ell}\to 0$. Thus,
the result \reff{s4f2a} follows since $z({\hat y_{\ell}}^*,x_k;\mu_k,\hat\rho_{\ell})\to\max\{A^Ty^*-c,0\}$ as $\hat\rho_{\ell}\to 0$.

In addition, since $s({\hat y_{\ell}}^*,x_k;\mu_k,\hat\rho_{\ell})-c+A^T{\hat y_{\ell}}^*>\mu_k$ for given $\mu_k>0$ and $\ell\in\{\ell|\hat\rho_{\ell+1}\le0.5\hat\rho_{\ell}\}$, and \bea
s({\hat y_{\ell}}^*,x_k;\mu_k,\hat\rho_{\ell})-c+A^T{\hat y_{\ell}}^*\to\max\{A^Ty^*-c,0\}~\hbox{as}~\hat\rho_{\ell}\to 0, \nn\eea
then $\max\{A^Ty^*-c,0\}\ge\mu_k>0$. That is, $y^*$ is infeasible to the problem \reff{prob1}, which by \cite{NocWri99,SunYua06,wright97,ye} implies that problem \reff{prob2} is unbounded. Noting that \reff{s4f2a} suggests that $y^*$ satisfies the stationary condition of the linear least square problem \bea
\min_y~\hf\|\max\{A^Ty-c,0\}\|^2, \nn\eea
$y^*$ is a point with the least $\ell_2$-norm of constraint violations of problem \reff{prob1}.  \eop

In the following analysis of this section, we suppose that the Step 1 of \refal{alg1} terminates finitely for every $k$. In order to analyze the convergence of \refal{alg1}, we also suppose that \refal{alg1} does not terminate finitely, and an infinite sequence $\{y_k\}$ is generated. 
Corresponding to the sequence $\{y_k\}$, we also have the sequence $\{\mu_k\}$ of barrier parameters, the sequence $\{\rho_k\}$ of penalty parameters, the sequence $\{x_k\}$ of the estimates of multipliers. In particular, $\{\mu_k\}$ is a monotonically decreasing sequence and tends to $0$, $\{\rho_k\}$ is a monotonically non-increasing sequence which either keeps unchanged after a finite number of steps or tends to $0$, \bea\dd\dd x_{k+1}=x_k+\frac{1}{\rho_{k+1}}(s(y_{k+1},x_{k};\mu_{k},\rho_{k+1})-c+A^Ty_{k+1}) \nn\\
\dd\dd=x_{k-1}+\frac{1}{\rho_{k}}(s(y_{k},x_{k-1};\mu_{k-1},\rho_{k})-c+A^Ty_{k})+\frac{1}{\rho_{k+1}}(s(y_{k+1},x_{k};\mu_{k},\rho_{k+1})-c+A^Ty_{k+1}) \nn\\
\dd\dd=x_0+\sum_{\ell=0}^k\frac{1}{\rho_{\ell+1}}(s(y_{\ell+1},x_{\ell};\mu_{\ell},\rho_{\ell+1})-c+A^Ty_{\ell+1}). \nn\eea
If the sequence $\{x_k\}$ is bounded, then $\frac{1}{\rho_{k+1}}(s(y_{k+1},x_{k};\mu_{k},\rho_{k+1})-c+A^Ty_{k+1})\to 0$ as $k\to\infty$, and $\{\rho_k\}$ is bounded away from zero.

\ble\label{lem41} If $\rho_k\to 0$, then any cluster point of $\{y_k\}$ is a Fritz-John point of problem \reff{prob1}. In particular, there exists an infinite subset ${\cal K}$ of indices such that for $k\in{\cal K}$ and $k\to\infty$, $y_k\to y^*$, $z_k\to z^*\ge 0$, $s_k\to s^*\ge 0$ and \bea s^*-c+A^Ty^*=0,\quad Az^*=0,\quad z^*\circ s^*=0, \label{s4f1}\eea
which shows that problem \reff{prob1} is feasible but problem \reff{prob2} is unbounded. 
\ele\prf Without loss of generality, we assume that $\{y_k\}$ is bounded. 
Because of the boundedness of $\{{\rho_k}{x_k}\}$, both $\{s_k\}$ and $\{z_k\}$ are bounded. Without loss of generality, we let $y_k\to y^*$, $z_k\to z^*$, $s_k\to s^*$ for $k\in{\cal K}$ and $k\to\infty$. Then $z^*\ge0$ and $s^*\ge 0$. Therefore, \reff{s4f1} follows immediately from \bea
\mu_k\to 0,~E_k^{primal}\le\mu_k,~E_k^{dual}\le\max\{\rho_{k},\mu_k\},~\hbox{and}~z_k\circ s_k={\rho_k}{\mu_k} e. \nn\eea
That is, $y^*$ is a Fritz-John point of problem \reff{prob1}. 

The equations in \reff{s4f1} show that, if $\rho_k\to 0$, \refal{alg1} will converge to a feasible point $y^*$ of \reff{prob1}. 
By the first part of the proof of \refl{le42n}, the finite termination of Step 1 implies that problem \reff{prob2} is strictly feasible. Thus, its set of solutions are unbounded since for any feasible point $x$ of problem \reff{prob2}, due to \reff{s4f1}, $x+\alpha z^*$ is feasible to problem \reff{prob2} and $c^T(x+\alpha z^*)=c^Tx$ for all $\alpha\ge 0$. \eop

In what follows, we prove the convergence of \refal{alg1} to a KKT point.
\ble\label{lem42} If $\rho_k$ is bounded away from zero, then $\{x_k\}$ is bounded, and every cluster point of $\{(y_k,x_k)\}$ is a KKT pair of problem \reff{prob1}. \ele\prf
Suppose that $\rho_{k}\ge\rho^*>0$ for all $k\ge 0$ and for some scalar $\rho^*$, then by Step 3 of \refal{alg1}, $\frac{1}{\|x_k\|_{\infty}}\ge\rho^*$. Thus, $\|x_k\|_{\infty}\le\frac{1}{\rho^*}$.

Since $\|x_k\|$ is bounded, $\lim_{k\to\infty}E_{k}^{dual}=0$. Thus, $E_{k}^{dual}\le\rho^*$ for all $k$ sufficiently large. Together with the facts $\mu_k\to 0$ and $$E_k^{primal}=\|\rho_kAx_{k}-\rho_k b\|_{(AA^T)^{-1}}\le\mu_k,$$
one has the result immediately. \eop

In summary, we have the following global convergence results on \refal{alg1}.
\bth One of following three cases will arise when implementing \refal{alg1}.

(1) The Step 1 does not terminate finitely for some $k\ge 0$, $\hat\rho_{\ell}\to 0$, either problem \reff{prob1} is unbounded, or problem \reff{prob1} is infeasible and problem \reff{prob2} is unbounded, and a point for minimizing the $\ell_2$ norm of constraint violations is found.

(2) The Step 1 terminate finitely for all $k\ge 0$, $\mu_k\to 0$ and $\rho_k\to 0$ as $k\to\infty$, problem \reff{prob2} is unbounded, problem \reff{prob1} is feasible and every cluster point of $\{y_k\}$ is a Fritz-John point of problem \reff{prob1}.

(3) The Step 1 terminate finitely for all $k\ge 0$, $\mu_k\to 0$ as $k\to\infty$, and $\rho_k$ is bounded away from zero, both problems \reff{prob1} and \reff{prob2} are feasible and every cluster point of $\{y_k\}$ is a KKT point of problem \reff{prob1}.
\eth\prf The results can be obtained straightforward from the preceding results Lemmas \ref{le41n}, \ref{le42n}, \ref{lem41}, and \ref{lem42} in this section. \eop

For reader's convenient, we summarize our global convergence results in Table \ref{tb1}.

\renewcommand\arraystretch{2}
\begin{table}[ht!b]
	\centering 
	\caption{The overview on the global convergence results of \refal{alg1}.}\label{tb1}
	\begin{tabular}{|c|c|c|>{\centering\arraybackslash}p{0.35\textwidth}|}	
		\hline 
		\multicolumn{1}{|c|}{\multirow{2}{*}{\refal{alg1}}}  &\multicolumn{3}{c|}{Results}\\
		\cline{2-4}
		&\multicolumn{1}{c|}{Dual LP \reff{prob1}} &\multicolumn{1}{c|}{Primal LP \reff{prob2}}  &\multicolumn{1}{c|}{Solution obtained} 
		\\
		\hline
		\multirow{2}{*}{$\hat\rho_{\ell}\to 0$, $\mu_k> 0$} 	
		& unbounded & - &-  \\	
		\cline{2-4}
		& infeasible & unbounded &\multirow{1}{0.35\textwidth}{A point for minimizing
			constraint violations of LP \reff{prob1}}\\			
		\hline	
		
		$\mu_k\to 0$, $\rho_k\to 0$
		& feasible
		& unbounded
		& A Fritz-John point of LP \reff{prob1}\\
		\hline	
		
		$\mu_k\to 0$, $\rho_k>0$
		& feasible
		& feasible
		& A KKT point\\
		\hline		
	\end{tabular}
\end{table}

\sect{Convergence rates and the complexity}

In this section, we concern about the convergence rate of \refal{alg1} under the situation that both problems \reff{prob1} and \reff{prob2} are feasible, which corresponds to the result (3) of the preceding global convergence theorem. Firstly, without any additional assumption, based on theory on convex optimization \cite{Nest18}, we prove that for given penalty parameter $\rho_k$, the convergence rate of the sequence of objective function values on the SBAL minimization subproblem is ${\large O}(\frac{1}{\ell})$, where $\ell>0$ is a positive integer which is also the number of iterations of the Step 1. Secondly, under the regular conditions on the solution, we show that the iterative sequence $\{{\hat y}_{\ell}\}$ on the SBAL minimization subproblem is globally linearly convergent. Finally, without loss of generality, by assuming that $\rho_k$ is small enough such that in Step 1, $E^{dual}_{k+1}\le\max\{\rho_k,\mu_k\}$ for given $\rho_k$, and using the preceding global linear convergence result, we can establish the iteration complexity of our algorithm. 
 
\bth For given $x_k$ and parameters $\mu_k$ and $\rho_k$, let $F_k(y)=L_B(y,x_k;\mu_k,\rho_k)$, $\{\hat y_{\ell}\}$ be a sequence generated by Step 1 of \refal{alg1} for minimizing $F_k(y)$, and $F_k^*=\inf_y F_k(y)$, $y_k^*=\hbox{\rm argmin}_y F_k(y)$. Then \bea
F_k(\hat y_{\ell})-F_k^*\le\frac{1}{2{\ell}}\|\hat y_0-y_k^*\|_{AA^T}^2, \eea
where $\hat y_0$ is an arbitrary starting point. \eth\prf It follows from \refl{lem22} that \bea
F_k(\hat y_{\ell+1})\dd\le\dd F_k(\hat y_{\ell})-\frac{1}{2}\|Az(\hat y_{\ell},x_k;\mu_k,\rho_k)-\rho_k b\|^2_{(AA^T)^{-1}} \nn\\
               \dd\le\dd F_k^*+\na F_k(\hat y_{\ell})^T(\hat y_{\ell}-y_k^*)-\frac{1}{2}\|Az(\hat y_{\ell},x_k;\mu_k,\rho_k)-\rho_k b\|^2_{(AA^T)^{-1}} \nn\\
               \dd=\dd F_k^*+\frac{1}{2}(\|\hat y_{\ell}-y_k^*\|_{AA^T}^2-\|\hat y_{\ell}-y_k^*-(AA^T)^{-1}\na F_k(\hat y_{\ell})\|_{AA^T}^2) \\
               \dd=\dd F_k^*+\frac{1}{2}(\|\hat y_{\ell}-y_k^*\|_{AA^T}^2-\|\hat y_{\ell+1}-y_k^*\|_{AA^T}^2), \nn\eea
where the second inequality follows from the convexity of $F_k(y)$, and the last equality is obtained by \reff{yiter}.
Thus, \bea
\sum_{t=1}^{\ell}(F_k(\hat y_{t})-F_k^*)\dd\le\dd\sum_{t=1}^{\ell}\frac{1}{2}(\|\hat y_{t-1}-y_k^*\|_{AA^T}^2-\|\hat y_{t}-y_k^*\|_{AA^T}^2) \nn\\
\dd=\dd\frac{1}{2}(\|\hat y_{0}-y_k^*\|_{AA^T}^2-\|\hat y_{\ell}-y_k^*\|_{AA^T}^2), \nn \eea
which implies $F_k(\hat y_{\ell})-F_k^*\le\frac{1}{2\ell}\|\hat y_{0}-y_k^*\|_{AA^T}^2$.
\eop

In order to derive the convergence rate of the iterative sequence $\{\hat y_{\ell}\}$ of our method for the subproblem, we need to prove some lemmas.
\ble\label{lem0614a} For given $x_k$ and parameters $\mu_k$ and $\rho_k$, let $F_k(y)=L_B(y,x_k;\mu_k,\rho_k)$.
Then for any $u,v\in\Re^m$, \bea (\na F_k(u)-\na F_k(v))^T(u-v)\ge\|\na F_k(u)-\na F_k(v)\|_{(AA^T)^{-1}}^2. \label{lem52f1}
\eea\ele\prf
For proving \reff{lem52f1}, we consider the auxiliary function \bea G_u(v)=F_k(v)-\na F_k(u)^Tv, \nn\eea
where $v$ is the variable and $u$ is any given vector. Then $\na G_u(u)=0$ and $\na^2 G_u(v)=\na^2 F_k(v)$, which means that $G_u(v)$ is convex as $F_k(v)$ and $u$ is precisely a global minimizer of $G_u(v)$. Therefore, we have a similar result to \refl{lem21} (1), that is, for every $w,v\in\Re^m$, \bea
G_u(w)\le G_u(v)+\na G_u(v)^T(w-v)+\hf(w-v)^TAA^T(w-v), \nn \eea
which implies $G_u(v)-G_u(u)\ge\frac{1}{2}\|\na G_u(v)\|_{(AA^T)^{-1}}^2$ for every $v\in\Re^m$. Because of $\na G_u(v)=\na F_k(v)-\na F_k(u)$, the preceding inequality is equivalent to \bea
F_k(v)-F_k(u)-\na F_k(u)^T(v-u)\ge\frac{1}{2}\|\na F_k(v)-\na F_k(u)\|_{(AA^T)^{-1}}^2. \label{lem52f2}\eea
Similarly, one can prove \bea
F_k(u)-F_k(v)-\na F_k(v)^T(u-v)\ge\frac{1}{2}\|\na F_k(v)-\na F_k(u)\|_{(AA^T)^{-1}}^2. \label{lem52f3}\eea
Summarizing two sides of \reff{lem52f2} and \reff{lem52f3} brings about our desired result.  \eop

In the subsequent analysis, let $y^*$ be the solution of problem \reff{prob1} and $x^*$ be the associated Lagrange multiplier vector, and $s^*=c-A^Ty^*$. Thus, $x^*\circ s^*=0$. We need the following blanket assumption.
\bas\label{ass2} Denote ${\cal I}=\{i=1,\ldots,n| x_i^*>0\}$. Suppose that the strict complementarity holds, and the columns of $A$ corresponding to the positive components of $x^*$ are linearly independent. That is, $x^*+s^*>0$ and $|{\cal I}|=m$, $B=A_{{\cal I}}A_{{\cal I}}^T$ is positive definite, where $|\cdot|$ is the cardinality of the set, $A_{{\cal I}}$ is a submatrix of $A$ consisting of $A_i,~i\in{\cal I}$. \eas

Under the \refa{ass2}, there exists a scalar $\delta>0$ such that, for $i\in{\cal I}$ and for all $\ell\ge 0$, $(s_{\ell i}+z_{\ell i})^{-1}z_{\ell i}\ge\delta>0$. Thus, for any $y\in\Re^m$, \bea
y^TA(S+Z)^{-1}ZA^Ty\dd\dd\ge y^T(A_{{\cal I}}(S_{{\cal I}}+Z_{{\cal I}})^{-1}Z_{{\cal I}}A_{{\cal I}}^T)y \nn\\
\dd\dd\ge\delta y^T(A_{{\cal I}}A_{{\cal I}}^T)y\ge\delta^{'}y^Ty\ge\delta^{''}y^TAA^Ty, \nn \eea
where $\delta^{'}\le\delta\lambda_{\min}(AA^T)$ and $\delta^{''}\le\frac{\delta^{'}}{\lambda_{\max}(AA^T)}<1$.
\ble For given $x_k$ and parameters $\mu_k$ and $\rho_k$, let $F_k(y)=L_B(y,x_k;\mu_k,\rho_k)$. Under the \refa{ass2}, there exists a scalar $\delta^{''}\in (0,1)$ such that, for any $u, v\in\Re^m$, \bea
\dd\dd (\na F_k(u)-\na F_k(v))^T(u-v) \nn\\
\dd\dd\ge\frac{1}{1+\delta^{''}}\|\na F_k(u)-\na F_k(v)\|_{(AA^T)^{-1}}^2+\frac{\delta^{''}}{1+\delta^{''}}\|u-v\|_{AA^T}^2. \eea
\ele\prf Let $G_k(y)=F_k(y)-\frac{1}{2}\delta^{''} y^TAA^Ty$. Then $G_k(y)$ and $\hf (1-\delta^{''})y^TAA^Ty-G_k(y)$ are convex, which suggests that $G_k(y)$ shares the similar properties with $F_k(y)$. Thus, the result of \refl{lem0614a} still holds for $G_k(y)$, i.e.,
for any $u,v\in\Re^m$, \bea (\na G_k(u)-\na G_k(v))^T(u-v)\ge\|\na G_k(u)-\na G_k(v)\|_{(AA^T)^{-1}}^2. \nn\eea
Due to $\na G_k(y)=\na F_k(y)-\delta^{''}AA^Ty$, the preceding inequality can be rewritten as \bea
\dd\dd (\na F_k(u)-\na F_k(v))^T(u-v) \nn\\
\dd\dd \ge\frac{1}{1-\delta^{''}}\|\na F_k(u)-\na F_k(v)-\delta^{''}AA^T(u-v)\|_{(AA^T)^{-1}}^2+\delta^{''}\|u-v\|_{AA^T}^2. \nn\eea
Thus, one has \bea \dd\dd (\na F_k(u)-\na F_k(v))^T(u-v) \nn\\
\dd\dd \ge\frac{1}{1+\delta^{''}}\|\na F_k(u)-\na F_k(v)\|_{(AA^T)^{-1}}^2+\frac{\delta^{''}}{1+\delta^{''}}\|u-v\|_{AA^T}^2,  \nn\eea
which completes our proof. \eop

Set $u=\hat y_{\ell}$ and $v=y_k^*$. Due to $\na F_k(y^*)=0$, \bea\na F_k(\hat y_{\ell})^T(\hat y_{\ell}-y_k^*)\ge\frac{1}{1+\delta^{''}}\|\na F_k(\hat y_{\ell})\|_{(AA^T)^{-1}}^2+\frac{\delta^{''}}{1+\delta^{''}}\|\hat y_{\ell}-y_k^*\|_{AA^T}^2. \nn\eea
The next result shows that sequence $\{\hat y_{\ell}\}$ can be of global linear convergence for the SBAL minimization subproblem.
\bth\label{th8n} Let $y_k^*=\hbox{\rm argmin} F_k(y)$. Under \refa{ass2}, there is a scalar $\tau\in (0,1)$ such that \bea
\|\hat y_{\ell}-y_k^*\|_{AA^T}^2\le\tau^{\ell}\|\hat y_{0}-y_k^*\|_{AA^T}^2. \nn\eea
That is, $\{\hat y_{\ell}\}$ is of global linear convergence to $y_k^*$. 
\eth\prf Note that
\bea \dd\dd\|\hat y_{\ell+1}-y_k^*\|_{AA^T}^2 \nn\\
\dd\dd=\|\hat y_{\ell}-(AA^T)^{-1}\na F_k(\hat y_{\ell})-y_k^*\|_{AA^T}^2 \nn\\
\dd\dd=\|\hat y_{\ell}-y_k^*\|_{AA^T}^2-{2}\na F_k(\hat y_{\ell})^T(\hat y_{\ell}-y_k^*)+\|\na F_k(\hat y_{\ell})\|_{(AA^T)^{-1}}^2 \nn\\
\dd\dd\le(1-\frac{2\delta^{''}}{1+\delta^{''}})\|\hat y_{\ell}-y_k^*\|_{AA^T}^2+(1-\frac{2}{1+\delta^{''}})\|\na F_k(\hat y_{\ell})\|_{(AA^T)^{-1}}^2 \nn\\
\dd\dd=\frac{1-\delta^{''}}{1+\delta{''}}\|\hat y_{\ell}-y_k^*\|_{AA^T}^2-\frac{1-\delta^{''}}{1+\delta^{''}}\|\na F_k(\hat y_{\ell})\|_{(AA^T)^{-1}}^2 \nn\\
\dd\dd\le\frac{1-\delta^{''}}{1+\delta{''}}\|\hat y_{\ell}-y_k^*\|_{AA^T}^2.  \nn\eea
By setting $\tau=\frac{1-\delta^{''}}{1+\delta{''}}$, the result follows immediately. \eop

Finally, based on the preceding global linear convergence result, we can obtain a new iteration complexity result on the algorithms for linear programs.
\bth Suppose that both problems \reff{prob1} and \reff{prob2} are feasible, and \refa{ass2} holds. For $\rho_0$ sufficiently small, if \refal{alg1} is terminated when $\mu_k<\epsilon$, where $\epsilon>0$ is a pre-given tolerance, then the iteration complexities of the MM methods for the subproblem and for problem \reff{prob1} are respectively \bea
T_{\rm MM}=O\left(\frac{1}{\ln\sqrt{\frac{\kappa_A+1}{\kappa_A-1}}}\ln \left(\frac{1}{\epsilon}\right)\right),\quad T_{\rm PDMM}=O\left(\frac{1}{\ln\sqrt{\frac{\kappa_A+1}{\kappa_A-1}}}\left(\ln \left(\frac{1}{\epsilon}\right)\right)^2\right).
\eea
\eth\prf Due to \refl{lemzp} (2), one has \bea \dd\dd\|Az(\hat y_{\ell+1},x_k;\mu_k,\rho_0)-\rho_0 b\| \nn\\
\dd\dd=\|\na_y L(\hat y_{\ell+1},x_k;\mu_k,\rho_0)-\na_y L(y_{k}^*,x_k;\mu_k,\rho_0)\| \nn\\
\dd\dd\le\|\hat y_{\ell+1}-y_k^*\|_{AA^T}. \nn \eea
In order to obtain $\|Az(\hat y_{\ell+1},x_k;\mu_k,\rho_0)-\rho_0 b\|\le\mu_k\le\epsilon$, by \reft{th8n}, $T_{\rm MM}$ should satisfy \bea
\sqrt{\tau}^{T_{\rm MM}}\|y_k-y_k^*\|_{AA^T}\le\epsilon, \nn\eea
where $y_k=\hat y_0$, $\tau$ is denoted in \reft{th8n} and can be replaced by $\tau=\frac{\kappa_A-1}{\kappa_A+1}$ ($\kappa_A=\lambda_{\max}(AA^T)/\lambda_{\min}(AA^T)$).
Thus, \bea T_{\rm MM}\ln\frac{1}{\sqrt{\tau}}\ge\ln\frac{\|y_k-y_k^*\|_{AA^T}}{\epsilon}. \nn\eea
That is, \bea T_{\rm MM}=O\left(\frac{1}{\ln\sqrt{\frac{\kappa_A+1}{\kappa_A-1}}}\ln\left(\frac{1}{\epsilon}\right)\right). \nn\eea

In addition, similarly, the number of iterations needed for driving $\mu_k<\epsilon$ is \bea
T_{\rm out}\ge\frac{1}{\ln\frac{1}{\gamma}}\ln\frac{\mu_0}{\epsilon}. \nn\eea Thus, we have the estimate on the total number of iterations \bea T_{\rm PDMM}=\sum_{k=1}^{T_{\rm out}}T_{\rm MM}=T_{\rm out}T_{\rm MM}=O\left(\frac{1}{\ln\sqrt{\frac{\kappa_A+1}{\kappa_A-1}}}\left(\ln \left(\frac{1}{\epsilon}\right)\right)^2\right), \eea
which completes our proof. \eop

\sect{Conclusion}

The simplex methods and the interior-point methods are two kinds of main and effective methods for solving linear programs. Relatively, the former is more inexpensive for every iteration but may require more iterations to find the solution, while the latter is more expensive for one iteration but the number of iterations may not be changed greatly with different problems. Theoretically, the iteration complexity of the simplex methods can be exponential on the sizes of linear programs, while the interior-point methods can be polynomial.

In this paper, we present a primal-dual majorization-minimization method for linear programs. The method is originated from the application of the Hestenes-Powell augmented Lagrangian method to the logarithmic-barrier problems. A novel barrier augmented Lagrangian (SBAL) function with second-order smoothness and strict convexity is proposed. Based the SBAL function, a majorization-minimization approach is introduced to solve the augmented Lagrangian subproblems. Distinct from the existing simplex methods and interior-point methods for linear programs, but similar to some alternate direction methods of multipliers (ADMM), the proposed method only depends on a factorization of the constant matrix independent of iterations which can be done in the preprocessing, and does not need any computation on step sizes, thus is much more inexpensive for iterations and can be expected to be particularly appropriate for large-scale linear programs. The global convergence is analyzed without prior assuming either primal or dual problem to be feasible. Under the regular conditions, based on theory on convex optimization, we prove that our method can be of globally linear convergence. The results show that the iteration complexity on our method is dependent on the conditioned number of the product matrix of the coefficient matrix and its transpose.

\


\end{document}